\documentclass[reqno,11pt]{amsart}

\newtheorem{theorem}{Theorem}[section]
\newtheorem{lemma}[theorem]{Lemma}

\theoremstyle{definition}

\theoremstyle{remark}
\newtheorem{remark}[theorem]{Remark}

\makeatletter
\def\dashint{\operatorname%
{\,\,\text{\bf--}\kern-.98em\DOTSI\intop\ilimits@\!\!}}
\makeatother

\newcommand\bR{\mathbb{R}}

 \newcommand{\mysection}[1]{\section{#1}
 \setcounter{equation}{0}}

\begin{document}

\title[On a paper by Dipierro, Savin, and Valdinoci]
{On the paper  
``All functions are locally s-harmonic up to a small error"
by Dipierro, Savin, and Valdinoci}
\author[N. V. Krylov]{N.V. Krylov}
\address[N. V. Krylov]{127 Vincent Hall, University of Minnesota,
 Minneapolis, MN, 55455}
\email{nkrylov@umn.edu}


\begin{abstract} 
We give an appropriate version of the result
in the paper by Dipierro, Savin, and Valdinoci
for different, not necessarily fractional, powers
 of the Laplacian.
\end{abstract}

\maketitle

\mysection{Introduction}

Let $\bR^{d}$ be a $d$-dimensional Euclidean space.
Denote $B_{r}=\{x\in\bR^{d}
:|x|<r\}$, and for integers $k=0,1,...$ let $C^{k}(\bar B_{r})$
be the set of $k$-times continuously differentiable
functions on $\bar B_{r}$ provided with the usual norm
$\|\cdot\|_{C^{k}(B_{r})}$
equal to the sum of maximal magnitudes of all
derivatives up to the $k$-th order.

Here is Theorem 1.1 of \cite{DSV_17}
in a slightly different form.

\begin{theorem}
                                        \label{theorem 10.10.1}
For any $k\in\{0,1,...\}$, $s\in(0,1)$,  $f\in C^{k}(\bar B_{r})$,
and $\varepsilon>0$, there exists $u\in C^{\infty}_{0}(\bR^{d})$
such that $(-\Delta)^{s}u=0$ in $B_{1}$ and
$$
\|f-u\|_{C^{k}(  B_{1})}\leq\varepsilon.
$$
\end{theorem}

Actually, Theorem \ref{theorem 10.10.1} may look like a generalization
of Theorem 1.1 of \cite{DSV_17} because the smoothness class of $u$'s
in the latter is much wider. However, by using mollifications
one easily deduces Theorem \ref{theorem 10.10.1}
from Theorem 1.1 of \cite{DSV_17}.

Many interesting consequences are extracted in \cite{DSV_17},
in particular, that $s$-harmonic functions in $B_{1}$
do not satisfy the Harnack inequality. 

One can rewrite Theorem \ref{theorem 10.10.1} in different terms.
 For real $n$ introduce
$$
K_{n}(x,y)=|x-y|^{n},\quad K_{n}g(y)=\int_{\bR^{d}}K_{n}(x-y)g(x)\,dx.
$$
Since the function $u$ from Theorem \ref{theorem 10.10.1}
has a representation $u=c_{d,s}K_{2s-d}g$,
where $g=(-\Delta)^{s}u$ and $c_{d,s}$ is an appropriate constant,
Theorem \ref{theorem 10.10.1} is equivalent to the
following.
\begin{theorem}
                                        \label{theorem 10.10.2}
For any $k\in\{0,1,...\}$, $s\in(0,1)$,  $f\in C^{k}(\bar B_{1})$,
and $\varepsilon>0$, there exists $g\in C^{\infty} (\bR^{d})$
having all bounded derivatives of any order
such that $g=0$ in $B_{1}$, $K_{2s-d}g\in  C^{\infty}_{0}(\bR^{d})$,
and
$$
\|f-K_{2s-d}g\|_{C^{k}(  B_{1})}\leq\varepsilon.
$$
\end{theorem}

The goal of this paper is to almost extend this result
from $n=2s-d$, $s\in(0,1)$, to arbitrary $n\in(-d,0)$
such that $n+d \not \in\{2,4,...\}$. By ``almost''
we mean that we will have $K_{n}g\in  C^{\infty} (\bR^{d})$
with each derivative of any order tending to zero
as $|x|\to\infty$ instead of
$K_{2s-d}g\in  C^{\infty}_{0}(\bR^{d})$, but on the other hand
we will have $g\in C^{\infty}_{0} (B_{4}\setminus\bar B_{3})$.
 Wider range of $n$ than 
$(-d,0)$ is also treated.

\mysection{Main result}

Fix a real number  $n$   
such that $n+d \not \in\{2,4,...\}$
and fix a $k\in\{0,1,...\}$. Define $S $
as the set of functions $f$ on $B_{1}$ each of which
has a representation on $B_{1}$ as $f=K_{n}g$,
where $g\in   C^{\infty}_{0}(B_{4}\setminus\bar B_{3})$.
Observe right away that if $n\in(-d,0)$,
then $K_{n}g$ is well defined in $\bR^{d}$
and is bounded along with each of its derivatives of any order,
which, in addition, tend to zero at infinity.
This justifies the last paragraph in the Introduction.
 
Here is our   result.
\begin{theorem}
                                          \label{theorem 9.10.1}
 The set $S$
is everywhere dense in $C^{k}(\bar B_{1})$.

\end{theorem}

To prove the theorem we need rather  rough estimates
from the following lemma proved in Section 
\ref{section 10.9.1}.
\begin{lemma}
                                          \label{lemma 9.10.1}
Take $r\in\{0,1,...\}$ such that $2r>k$. Then
there exist a constant $N=N(d,k,r)$ such that,
for any $x\in B_{4}\setminus\bar B_{3}$ and real $m$,
  we have
\begin{equation}
                                                \label{9.10.1}
\|K_{m}(x,\cdot)\|_{C^{ k}(B_{1})}
\leq N^{r}(m^{2r}+N^{r}).
\end{equation}

\end{lemma}

{\bf Proof of Theorem \ref{theorem 9.10.1}}. Let $\bar S$ 
be the closure of 
$S$ in $C^{k}(\bar B_{1})$.
Then $K(x,\cdot)\in \bar S$ for any $x\in B_{4}\setminus\bar  B_{3}$,
which is seen from the formula
$$
K_{n}(x,y)=\lim_{m\to\infty}\int_{\bR^{d}}K_{n}(x+z/n,y)\zeta 
(z )\,dz
$$
$$
=\lim_{m\to\infty}\int_{B_{4}\setminus\bar  B_{3}}
K_{n}(z,y)\zeta_{m}
(z-x)\,dz,
$$
where $\zeta_{m}(z)=m^{d}\zeta(mz)$ and $\zeta$ is any
$C^{\infty}_{0}(\bR^{d})$ function which integrates to one.

Now, since partial derivatives are limits of finite differences,
$\Delta^{m}_{x}K_{n}(x,\cdot)\in \bar S$ for any $m=0,1,...$
and $x\in B_{4}\setminus \bar B_{3}$,
where $\Delta_{x}$ is the Laplacian applied with respect to
the $x$ variable. Observe that $K_{n}(x-y)$ is a radial function
if the origin is taken to be $y$. We know
how the Laplacian   looks on radial functions,
so that, for $m\geq1$,
\begin{equation}
                                                \label{9.10.2}
\Delta _{x}K_{n} =n(n-1)  K_{n-2}+(d-1)n
 K_{n-2}=
n(n+d-2) K_{n-2}.
\end{equation}
It follows  that
$\Delta^{m}_{x}K_{n}(x,\cdot)=c_{m}K_{n-2m}(x,\cdot)$,
where $c_{m}$ are some constants,
and our   assumption
that $n+d \not \in\{2,4,...\}$ implies that $c_{m}\ne0$.
Thus
$$
K_{n-2m}(x,\cdot)\in \bar S\quad m=0,1,....
$$

Next, owing to Lemma \ref{lemma 9.10.1},
 for any $x\in B_{4}\setminus \bar B_{3}$, $r\in\{0,1,...\}$,
and real $t$,
 the series
$$
\sum_{m=0}^{\infty}\frac{1}{m!}t^{m}K_{n-2r-2m}(x,\cdot)
$$
converges in $C^{k}(\bar B_{1})$. Thus, for any $t$ 
and $x\in B_{4}\setminus \bar B_{3}$,
$$
K_{n-2r}(x-\cdot)\exp(-t/|x-\cdot|^{2})\in \bar S.
$$
Obviously,
the integral
$$
\int_{0}^{\infty}t^{m}K_{n-2r}(x-\cdot)\exp(-t/|x-\cdot|^{2})\,dt
$$
converges in $C^{k}(\bar B_{1})$ for any
$x\in B_{4}\setminus \bar B_{3}$ and $m>-1$. Thus,
$$
K_{n-2r+2m-2}(x-\cdot)\in \bar S.
$$
For any real $j $ it is possible to find
(many) $r\in\{0,1,...\}$ and real $m>-1$ such that
$n-2r+2m-2=j$. Hence, for any
$x\in B_{4}\setminus \bar B_{3}$
$$
K_{j}(x-\cdot)\in \bar S\quad \forall j\in\bR.
$$
For $j\in\{0,1,...\}$ the functions 
$K_{j}(x-y)$ are polynomials with respect to
$x$ with the coefficients depending on $y$.
These coefficients can be found by differentiating
$K_{j}(x-y)$ with respect to $x\in B_{4}\setminus\bar B_{3}$
 and it follows that
they belong to $\bar S$. But then, for any $x$, 
$K_{j}(x-\cdot)$ is in $\bar S$  being just a linear
combination of elements of $\bar S$. Thus,
$$
K_{j}(x-\cdot)\in \bar S\quad \forall j=0,1,..., x\in\bR^{d}.
$$

Furthermore, clearly
the series of polynomials
$$
\sum_{m=0}^{\infty}\frac{1}{m!}t^{-m}K_{2m}(x-\cdot)
$$
converges in $C^{k}(\bar B_{1})$ for any
$x\in \bR^{d}$ and $t>0$. It follows
that
$$
\exp(-|x-\cdot|^{2}/t)\in \bar S.
$$

The rest is routine. Take an $f\in C^{k}(\bar B_{1})$,
continue it to the whole $\bR^{d}$ preserving its class
and making it vanish outside $B_{2}$,
and for $t>0$ define
$$
f_{t}(y)= t^{-d/2}\int_{\bR^{d}}
f(x)\exp(-|x-y|^{2}/t)\,dx.
$$
Obviously the integral converges in $C^{k}(\bar B_{1})$,
so that $f_{t}\in \bar S$.
We know that, if a constant $c$   is chosen appropriately,
then $cf_{t}\to f$ in $C^{k}(\bar B_{1})$ as $t\downarrow 0$.
Since $cf_{t}\in \bar S$, $f \in \bar S$
and the theorem is proved. 

\mysection{Proof of Lemma \protect\ref{lemma 9.10.1}}

                                        \label{section 10.9.1}

Take $p\in(1,\infty)$ and $u\in W^{2}_{p}(B_{2})$.
By representing it as $h+v$, where $h$ is harmonic
and $v=0$ on $\partial B_{2}$ it is easy to conclude that
$\|u\|_{W^{2}_{p}(B_{1})}\leq N(\|u\|_{L_{p}(B_{2})}
+\|\Delta u\|_{L_{p}(B_{2})})$, where and below by $N$ we denote
generic constants $\geq1$ depending only
on $d$, $p$, and $r$ (which is coming later). This fact also holds
for $B_{\rho_{1}}\subset \bar B_{\rho_{1}}\subset   B_{\rho_{2}}$
in place of $B_{1}$ and $B_{2}$ with a constant, this time,
also depending on $\rho_{1}$ and $\rho_{2}$.
By   applying these estimates to the derivatives of $u$, we see that,
for any $r=0,1,...$,
$$
\|u\|_{W^{2r}_{p}(B_{1})}
\leq N\sum_{i=0}^{r}\|\Delta^{i} u\|_{L_{p}(B_{2})}.
$$

  Then
 \eqref{9.10.2} and the estimate $|n(n+d-2)|\leq N(n^{2}+1)$ 
  show  that, for $|y|\leq2$ and $3\leq|x|\leq4$,
$$
|\Delta^{i} K_{m}(x,y)|\leq c^{i}(m^{2}+1)((m-2)^{2}+1)
\cdot...\cdot((m-2i+2)^{2}+1)K_{m-2i}(x,y)
$$
$$
\leq c^{i}(m^{2}+1)((m-2)^{2}+1)
\cdot...\cdot((m-2i+2)^{2}+1).
$$
Since, for $j=0,...,r$, we have $(m-2j)^{2}+1\leq 2m^{2}+N$,
it holds  that
$$
\|K_{m}(x,\cdot)\|_{W^{2r}_{p}(B_{1})}
\leq N^{r}(2m^{2}+N)^{r}.
$$

Finally, by embedding theorems
$W^{2r}_{p}(B_{1})\subset C^{k}(\bar B_{1})$ if
$(2r-k)p>d$, which holds indeed for large enough $p$.
This proves the lemma.

We conclude the paper with the following.
\begin{remark}
In \cite{DSV_17} the authors prove that 
the set of polynomials is dense 
in $C^{k}(\bar B_{1})$. Actually, this fact follows from
Lemma 3 and Theorem 3 in Section 3 of
an old paper \cite{Gr_41}. By using Bernstein's
polynomials this fact is also follows from
Section 3.8 of \cite{Kr_95}.

\end{remark}

\end{document}